\documentclass{article}

\newcommand{\textlineskip}{\baselineskip=11.5dd plus0.5dd minus0.2dd}

\def\qed{\hbox{${\vcenter{\vbox{
   \hrule height 0.4pt\hbox{\vrule width 0.4pt height 6pt
   \kern5pt\vrule width 0.4pt}\hrule height 0.4pt}}}$}}

\catcode`\@=11

\font\twrm=cmr12

\font\eightrm=cmr8

\newcounter{subsect}

\def\section#1{\par \bigskip \setcounter{subsect}{0}
\setcounter{equation}{0}
\addtocounter{section}{1}\begin{center}\thesection.
\uppercase{#1} \end{center} \par \smallskip }

\def\subsection#1{\par \medskip
\addtocounter{subsect}{1}\begin{center}{ \eightrm
\thesection.\thesubsect. \uppercase{#1}} \end{center} \par
\smallskip }

\usepackage{amsmath}
\usepackage{amsfonts}
\newcommand{\abs}[1]{{\left| {#1} \right|}}
\newcommand{\R}{\mathbb{R}}

\begin{document}

\normalsize\textlineskip
\thispagestyle{empty}

 \begin{center} {\twrm
%
WEIGHTED BANACH SPACES\\
\vspace*{0.035truein}OF HOLOMORPHIC FUNCTIONS WITH\\
\vspace*{0.035truein}LOG-CONCAVE WEIGHT FUNCTION
}\\
\vspace*{2cc}
{\eightrm Martin~At.~Stanev\\University of Forestry, Sofia}
\end{center}

\vspace*{30dd}

{\parindent0pt \footnotesize \leftskip20pt
\rightskip20pt  \baselineskip10pt
Some theorems on convex functions are proved and an application of these theorems in the theory of 
weighted Banach spaces of holomorphic functions is investigated, too. We prove that $H_v(G)$ and $H_{v_0}(G)$ are exactly the same spaces as $H_w(G)$ and $H_{w_0}(G)$ where $w$ is the smallest log-concave majorant of $v$. This investigation is based on the theory of convex functions and some specific properties of the weighted banach spaces of holomorphic functions under considaration. 
\smallskip

{\bf  2010
Math.\ Subject Classification} 46E15, 46B04.
\par
}

\vspace*{16dd}

\section{Introduction}
Let $\mathbb C$ be the complex plane and 
\[
G=\{ z=x+iy|\; x\in(-\infty;+\infty), y\in(0;+\infty)\}\subset\mathbb C
\] 
be the upper half plane of $\mathbb C$. The function 
$v: G\to(0;+\infty)$ is such that $v(z)=v(x+iy)=v(iy)$, 
$\forall z=x+iy\in G$, and
\[
\inf_{y,y\in[\frac1c,c]}v(iy)>0,\qquad \forall c>1.\eqno(1.1)
\]
We define
\[
\varphi_v(y)=(-1)\ln v(iy),\quad y\in(0;+\infty).
\]

Thus the property~(1.1) is reformulated as the folowing property of $\varphi_v(y)$
\[
\sup_{y,y\in[\frac1c,c]}\varphi_v(y)<+\infty,\qquad \forall c>1.\eqno(1.1')
\]

The weighted Banach spaces of holomorphic functions $H_v(G)$ and $H_{v_0}(G)$ are defined as follows
\begin{itemize}
\item $f\in H_v(G)$ iff $f$ is holomorphic on $G$ and is such that
\[
\parallel f\parallel_v=\sup_{z,z\in \mathcal G}v(z)\abs{f(z)},
\]
\item $f\in H_{v_0}(G)$ iff $f\in H_v(G)$ and $f$ is such that $\forall \varepsilon>0$ there exists a compact $\mathcal K_{\varepsilon}\subset G$ for which
\[
\sup_{z\in G\setminus \mathcal K_{\varepsilon}} v(z)\abs{f(z)}<\varepsilon.
\]
\end{itemize}
Thus, we use notations used in [1,2,3,4,5]. 

In [1], [2] authors find the isomorphic classification of the spaces $H_v(G)$ and $H_{v_0}(G)$ when the weight function $v$ meets some growth conditions. 

In [3], [4] are studied weighted composition operators between weighted spaces of holomorphic functions on the unit disk of the complex plane and the associated weights are used in order to estimate the norm of the weighted composition operators. 

In [5] are studied the associated weights.

This paper is about the weights that have some of the properties of the associated weights. We prove that $H_v(G)$ and $H_{v_0}(G)$ are exactly the same spaces as $H_w(G)$ and $H_{w_0}(G)$ where $w$ is the smallest log-concave majorant of $v$. Here,  the smallest log-concave majorant of $v$ is exactly the associated weight but in case of other weighted spaces this coincidation might not take place. Our work is based on the theory of convex functions and some specific properties of the weighted banach spaces of holomorphic functions under consideration.

The results of this paper are communicated on the conferences [7] and [8].

\section{Definitions and notations}
Let $\Phi$ be the set of functions $\varphi$ such that $\varphi\in\Phi$ iff the following conditions are fullfield
\begin{itemize}
\item $\varphi:(0;+\infty)\to\R$ and
\item there exists a real number $a$ such that
\[
\inf_{x\in(0;+\infty)}\; \bigl(\varphi(x)-ax\bigr)>-\infty.
\]
\end{itemize}
Note that $-\infty<\varphi(x)<+\infty$,  $\forall x\in(0;+\infty)$ and $\forall \varphi\in\Phi$.

We denote by $\widehat{a}_{\varphi}$ the limit inferior
\[
\widehat{a}_{\varphi}=\liminf_{x\to+\infty}\; \frac{\varphi(x)}{x}, \quad \varphi\in\Phi.
\]

If $\varphi\in\Phi$ then 
\begin{itemize}
\item $\widehat{a}_{\varphi}\in\R\cup\{+\infty\}$, $\widehat{a}_{\varphi}>-\infty$,
\item $\widehat{a}_{\varphi}=\sup\bigl\{ a\bigl|\; a\in\R,\; \inf\limits_{x\in(0;+\infty)}\; (\varphi(x)-ax)>-\infty \bigr.\bigr\} $
\end{itemize}

If $\varphi\in\Phi$ is convex on $(0;+\infty)$ then
\[
\widehat{a}_{\varphi}=\lim\limits_{x\to+\infty}\; \frac{\varphi(x)}{x}
\]

Let $\Phi_1$, $\Phi_2$, $\Phi_3$ be the following subsets of $\Phi$
\begin{align*}
\Phi_1&=\bigl\{\varphi\bigl|\; \varphi\in\Phi,\; \widehat{a}_{\varphi}=+\infty \bigr.\bigr\}\\
\Phi_2&=\bigl\{\varphi\bigl|\; \varphi\in\Phi,\; \widehat{a}_{\varphi}<+\infty, \;\liminf_{x\to+\infty} (\varphi(x)-\widehat{a}_{\varphi}x)=-\infty \bigr.\bigr\}\\
\Phi_3&=\bigl\{\varphi\bigl|\; \varphi\in\Phi,\; \widehat{a}_{\varphi}<+\infty, \;\liminf_{x\to+\infty} (\varphi(x)-\widehat{a}_{\varphi}x)>-\infty \bigr.\bigr\}
\end{align*}

Note that $\Phi_1$, $\Phi_2$, $\Phi_3$ are mutually disjoint sets and
$\Phi_1\cup\Phi_2\cup\Phi_3=\Phi$. 

If $\varphi\in\Phi_2\cup\Phi_3$ is convex on $(0;+\infty)$ then
\[
\liminf\limits_{x\to+\infty} (\varphi(x)-\widehat{a}_{\varphi}x)=\lim\limits_{x\to+\infty} (\varphi(x)-\widehat{a}_{\varphi}x).
\]

Note that a function $\varphi\in\Phi$ is not necessarily continuous. Here,  the function $\varphi$ that belongs to $\Phi$, does not have to meet any conditions beside those of the definition of $\Phi$, $\Phi_1$, $\Phi_2$, $\Phi_3$. 
There are a number of simple functions that belong to $\Phi$, $\Phi_1$, $\Phi_2$, $\Phi_3$ and for example
\begin{itemize}
\item if $\varphi_1(x)=x^2$,  $\forall x\in(0;+\infty)$, then $\varphi_1\in\Phi_1$;
\item if $\varphi_2(x)=x-\sqrt{x}$, $\forall x\in(0;+\infty)$, then $\varphi_2\in\Phi_2$;
\item if $\varphi_3(x)=x^{-1}$, $\forall x\in(0;+\infty)$, then $\varphi_3\in\Phi_3$;
\end{itemize}
and the functions $\varphi_1(x)$, $\varphi_2(x)$, $\varphi_3(x)$ are all convex on $(0;+\infty)$.

Let $\varphi\in\Phi$ and let 
\[
M_{\varphi}=\bigl\{(a,b)\bigl|\; a\in\R, b\in\R, \; \inf\limits_{t\in(0;+\infty)}\; (\varphi(t)-at)>b \bigr.\bigr\}.
\]
The function $\varphi^{**}:(0;+\infty)\to\R$ is defined as
\[
\varphi^{**}(x)=\sup_{(a,b)\in M_{\varphi}}\; (ax+b).
\]
Thus, $\varphi^{**}$ is the second  Young-Fenhel conjugate of $\varphi$
and it is the biggest convex minorant of $\varphi$.

\section{Main Results}
{\bf Theorem~3.1}
Let $\varphi\in\Phi$ and $\psi\in\Phi$. If $\psi$ is convex on $(0;+\infty)$ then 
\[
\inf_{x\in(0;+\infty)}\; \bigl(\varphi(x)-\psi(x)\bigr)=
\inf_{x\in(0;+\infty)}\; \bigl(\varphi^{**}(x)-\psi(x)\bigr)
\]

{\bf Theorem~3.2}
Let $\varphi\in\Phi$ and $\psi\in\Phi$. If $\psi$ is convex on $(0;+\infty)$   and the right-sided limit $\lim\limits_{x\to0^{+}} \bigl(\varphi(x)-\psi(x) \bigr)=+\infty$ then
\[
\lim_{x\to0^{+}} \bigl(\varphi^{**}(x)-\psi(x) \bigr)=+\infty
\]

{\bf Theorem~3.3}
Let $\varphi\in\Phi$, $\psi\in\Phi\setminus\Phi_3$. 
If $\psi$ is convex on $(0;+\infty)$ and the limit value $\lim\limits_{x\to+\infty} \bigl(\varphi(x)-\psi(x) \bigr)=+\infty$ then
\[
\lim_{x\to+\infty} \bigl(\varphi^{**}(x)-\psi(x) \bigr)=+\infty
\]

{\bf Example~3.1}
 Theorem~3.1 does not hold with the functions
\[
\begin{aligned}
\varphi(x)&=\min\{x,1\}+1,&&\nonumber\\
\psi(x)&=\frac{x}{x+1},&&\forall x\in(0;+\infty)\nonumber
\end{aligned}
\]
Note that $\varphi\in\Phi$, $\psi\in\Phi$, the function $\psi$ is {\it not} convex on $(0;+\infty)$, $\varphi^{**}(x)=1$,  $\forall x\in(0;+\infty)$, and  
\[
1=\inf\limits_{x\in(0;+\infty)}\; \bigl(\varphi(x)-\psi(x)\bigr)\neq
\inf\limits_{x\in(0;+\infty)}\; \bigl(\varphi^{**}(x)-\psi(x)\bigr)=0.
\]

{\bf Example~3.2}
 Theorem~3.2 does not hold with the functions 
\begin{align*}
\varphi(x)&=\frac{1}{x^2}+\frac{1}{x}\sin\frac{1}{x}+\frac{2}{x},\\
\psi(x)&=\varphi(x)-\frac{2}{x},
\end{align*}
where $x\in(0;+\infty)$. Note that $\varphi\in\Phi$, $\psi\in\Phi$, the function $\psi$ is {\it not} convex on $(0;+\infty)$ and 
\[
+\infty=\lim\limits_{x\to0^{+}} \bigl(\varphi(x)-\psi(x) \bigr)>\liminf_{x\to0^{+}} \bigl(\varphi^{**}(x)-\psi(x) \bigr)
\]
This is proved in the Proposition~4.1.

{\bf Example~3.3}
Theorem~3.3 does not hold with the functions  
\begin{align*}
\varphi(x)=x^2+x\sin x +2x, \\
\psi(x)=\varphi(x)-2x,
\end{align*}
where $x\in(0;+\infty)$. Note that $\varphi\in\Phi$, $\psi\in\Phi\setminus\Phi_3$, the function $\psi$ is {\it not} convex on $(0;+\infty)$ and
\[
+\infty=\lim\limits_{x\to+\infty} \bigl(\varphi(x)-\psi(x) \bigr)>\liminf_{x\to+\infty} \bigl(\varphi^{**}(x)-\psi(x) \bigr)
\]
This is proved in the Proposition~4.2.

{\bf Corollary~3.1}
Let $\varphi\in\Phi$ and $\psi\in\Phi$. If $\varphi$ and $\psi$ are such that the right-sided limit 
\[
\lim\limits_{x\to0^{+}} \bigl(\varphi(x)-\psi(x) \bigr)=+\infty,
\]
then
\[
\lim_{x\to0^{+}} \bigl(\varphi^{**}(x)-\psi^{**}(x) \bigr)=+\infty\eqno(3.1)
\]

{\it Proof.} Note that $\psi^{**}\le\psi$ and 
\[
\lim\limits_{x\to0^{+}} \bigl(\varphi(x)-\psi^{**}(x) \bigr)\ge \lim\limits_{x\to0^{+}} \bigl(\varphi(x)-\psi(x) \bigr)=+\infty.
\]
The Theorem~3.2 is applied to $\varphi$, $\psi^{**}$ and thus the limit~(3.1) is proved. \qed

{\bf Corollary~3.2}
Let $\varphi\in\Phi$, $\psi\in\Phi\setminus\Phi_3$. 
If $\varphi$ and $\psi$ are such that  
\[
\lim\limits_{x\to+\infty} \bigl(\varphi(x)-\psi(x) \bigr)=+\infty,
\]
then
\[
\lim_{x\to+\infty} \bigl(\varphi^{**}(x)-\psi^{**}(x) \bigr)=+\infty\eqno(3.2)
\]

{\it Proof.} Note that
\begin{itemize}
\item $\psi^{**}\in\Phi\setminus\Phi_3$ by the Lemma~4.1;
\item $\psi^{**}\le\psi$.
\end{itemize}
The Theorem~3.3 is applied to $\varphi$, $\psi^{**}$ and thus the limit~(3.2) is proved. \qed

{\bf Example~3.3}
Let
\begin{align*}
\varphi(x)&=x^2+x,\\
\psi(x)&=\begin{cases}
3x-1,&x\in(0,1]\\
5-3x,&x\in(1,2]\\
x^2+x-7,&x\in(2,+\infty)
\end{cases}
\end{align*}
So, $\varphi\in\Phi$, $\psi\in\Phi$ and
\begin{itemize}
\item $\varphi$ is convex on $(0,+\infty)$, and therefore $\varphi^{**}=\varphi$,
\item $\psi$ is {\it not} convex on $(0,+\infty)$ and
\[
\psi^{**}(x)=\begin{cases}
-1,&x\in(0,2]\\
x^2+x-7,&x\in(2,+\infty)
\end{cases}
\]
\end{itemize}
A direct calculation shows that
\[
\inf_{x\in(0;+\infty)}\; \bigl(\varphi(x)-\psi(x)\bigr)=0\neq1=\inf_{x\in(0;+\infty)}\; \bigl(\varphi^{**}(x)-\psi^{**}(x)\bigr).
\]
Thus there is no any analog of the Theorem~3.1 involving $\varphi^{**}$ and $\psi$ in such a way as in Corollaries~3.1 and~3.2. \qed

\section{Auxiliary results}
{\bf Proposition~4.1} Let 
\begin{align*}
\varphi(x)&=\frac{1}{x^2}+\frac{1}{x}\sin\frac{1}{x}+\frac{2}{x},\\
\psi(x)&=\varphi(x)-\frac{2}{x},
\end{align*}
where $x\in(0;+\infty)$. The functions $\varphi$ and $\psi$ are such that $\varphi\in\Phi$, $\psi\in\Phi$, the function $\psi$ is not convex on $(0;+\infty)$ and 
\[
+\infty=\lim\limits_{x\to0^{+}} \bigl(\varphi(x)-\psi(x) \bigr)>\liminf_{x\to0^{+}} \bigl(\varphi^{**}(x)-\psi(x) \bigr)
\]

{\it Proof.} 
The function $\psi$ is such that
\[
\psi(x)\ge\begin{cases}
\frac{1}{x^2}-\frac{1}{x},& x\in(0;1)\\
\frac{1}{x^2},&x\in[1;+\infty)
\end{cases}
\]
so, $\psi(x)\ge0$ $\forall x\in(0;+\infty)$ and this implies that $\psi\in\Phi$. 

The function $\varphi\in\Phi$ because of the inequality $\varphi\ge\psi$.

Note that the limit value
\[
\lim\limits_{x\to0^{+}} \bigl(\varphi(x)-\psi(x) \bigr)=\lim\limits_{x\to0^{+}} \frac{2}{x}=+\infty.
\]

Let
\[
x_k=\frac{1}{\frac{3\pi}{2}+2k\pi},\quad 
\widetilde{x}_k=\frac{1}{\frac{5\pi}{2}+2k\pi}
\]
where $k=0,1,2,\dots$. Note that $x_k>\widetilde{x}_k>x_{k+1}>0$, $\lim\limits_{k\to+\infty} x_k=0$, and the harmonic mean of $x_k$, $x_{k+1}$ is equal to $\widetilde{x}_{k}$. 

A direct computation shows that the second derivative
$\psi''(\widetilde{x}_0)<0$. Therefore $\psi$ is not convex on $(0;+\infty)$.

Let 
\[
f(x)=\frac{1}{x^2}+\frac{1}{x},\quad \forall x\in(0;+\infty).
\]
Note that the function $f$ is convex on $(0;+\infty)$ and $f(x)\le\varphi(x)$, $\forall x\in(0;+\infty)$. So, $f$ is a convex minorant of $\varphi$ and thus $f\le\varphi^{**}$.

Therefore, $f(x_k)\le\varphi^{**}(x_k)\le\varphi(x_k)=f(x_k)$ and this implies that  
\[
f(x_k)=\varphi^{**}(x_k), \quad \forall k=1,2,3,\dots.
\]
Furhtermore, 
\[
\psi(\widetilde{x}_k)=f(\widetilde{x}_k)\le\varphi^{**}(\widetilde{x}_k)\le \frac{x_k-\widetilde{x}_k}{x_k-x_{k+1}}\varphi^{**}(x_{k+1})+\frac{\widetilde{x}_k-x_{k+1}}{x_k-x_{k+1}}\varphi^{**}(x_k)
\]
because of the convexity of $\varphi^{**}$. Thus
\[
0\le\varphi^{**}(\widetilde{x}_k)-\psi(\widetilde{x}_k)\le \frac{x_k-\widetilde{x}_k}{x_k-x_{k+1}}f(x_{k+1})+\frac{\widetilde{x}_k-x_{k+1}}{x_k-x_{k+1}}f(x_k)-f(\widetilde{x}_k)
\]
Note that after some simple computations we obtain
\[
\frac{x_k-\widetilde{x}_k}{x_k-x_{k+1}}f(x_{k+1})+\frac{\widetilde{x}_k-x_{k+1}}{x_k-x_{k+1}}f(x_k) - f(\widetilde{x}_k)=(3+\widetilde{x}_k)\pi^2
\]

Consequently $0\le\varphi^{**}(\widetilde{x}_k)-\psi(\widetilde{x}_k)\le (3+\widetilde{x}_k)\pi^2$, $\forall k=1,2,3,\dots$ and
\[
\liminf_{x\to0^{+}} \bigl(\varphi^{**}(x)-\psi(x) \bigr)<+\infty.\qquad\qed
\]

{\bf Proposition~4.1} Let
\begin{align*}
\varphi(x)=x^2+x\sin x +2x, \\
\psi(x)=\varphi(x)-2x,
\end{align*}
where $x\in(0;+\infty)$. The functions $\varphi$ and $\psi$ are such that $\varphi\in\Phi$, $\psi\in\Phi\setminus\Phi_3$, the function $\psi$ is not convex on $(0;+\infty)$ and
\[
+\infty=\lim\limits_{x\to+\infty} \bigl(\varphi(x)-\psi(x) \bigr)>\liminf_{x\to+\infty} \bigl(\varphi^{**}(x)-\psi(x) \bigr)
\]

{\it Proof.} The function $\psi$ is such that 
\[
\psi(x)\ge\begin{cases}
x^2,& x\in(0;\pi)\\
x^2-x,&x\in[\pi;+\infty)
\end{cases}
\]
Therefore $\psi\ge0$ and thus $\psi\in\Phi$.

The limit value
\[
\widehat{a}_{\psi}=\liminf_{x\to+\infty} \;\frac{\psi(x)}{x}\ge\liminf_{x\to+\infty} \;\frac{x^2-x}{x}=+\infty
\] 
Therefore $\widehat{a}_{\psi}=+\infty$ and thus $\psi\in\Phi_1\subset\Phi\setminus\Phi_3$.

The function $\varphi\in\Phi$ because of both $\varphi\ge\psi$ and $\psi\in\Phi$.

The limit value 
\[
\lim\limits_{x\to+\infty}\bigl(\varphi(x)-\psi(x) \bigr)=\lim\limits_{x\to+\infty} 2x=+\infty.
\]

Let
\[
x_k=\frac{3\pi}{2}+2k\pi,\quad 
\widetilde{x}_k=\frac{5\pi}{2}+2k\pi
\]
where $k=0,1,2,\dots$. Note that, if $k=1,2,3,\dots$ then $0<x_k<\widetilde{x}_k<x_{k+1}$, $x_k+x_{k+1}=2\widetilde{x}_{k}$ and $\lim\limits_{k\to+\infty} x_k=+\infty$.

A direct computation shows that the second derivative $\psi''(\widetilde{x}_0)<0$
and so the function $\psi$ is not convex on $(0;+\infty)$.

Let 
\[
f(x)=x^2+x,\quad \forall x\in(0;+\infty).
\]
The function $f$ is convex on $(0;+\infty)$ and $f\le\varphi$.
So,  $f$ is convex minorant of $\varphi$ and thus 
$f\le\varphi^{**}$.

Therefore $f(x_k)\le\varphi^{**}(x_k)\le\varphi(x_k)=f(x_k)$ and this implies 
\[
f(x_k)=\varphi^{**}(x_k), \quad \forall k=1,2,3,\dots.
\]
Furthermore,
\[
\psi(\widetilde{x}_k)=f(\widetilde{x}_k)\le\varphi^{**}(\widetilde{x}_k)\le \frac{x_{k+1}-\widetilde{x}_k}{x_{k+1}-x_{k}}\varphi^{**}(x_k)+\frac{\widetilde{x}_k-x_k}{x_{k+1}-x_k}\varphi^{**}(x_{k+1})
\]
because of the convexity of $\varphi^{**}$. Thus
\[
0\le\varphi^{**}(\widetilde{x}_k)-\psi(\widetilde{x}_k)\le \frac{x_{k+1}-\widetilde{x}_k}{x_{k+1}-x_{k}}f(x_k)+\frac{\widetilde{x}_k-x_k}{x_{k+1}-x_k}f(x_{k+1})-f(\widetilde{x}_k)
\]

Note that after some simple computations we obtain
\[
\frac{x_{k+1}-\widetilde{x}_k}{x_{k+1}-x_{k}}f(x_k)+\frac{\widetilde{x}_k-x_k}{x_{k+1}-x_k}f(x_{k+1})-f(\widetilde{x}_k)=\pi^2, \quad \forall k=1,2,3,\dots.
\]

Consequently, $0\le\varphi^{**}(\widetilde{x}_k)-\psi(\widetilde{x}_k)\le \pi^2$, $\forall k=1,2,3,\dots$  and
\[
\liminf\limits_{x\to+\infty} \bigl(\varphi^{**}(x)-\psi(x) \bigr)<+\infty.\qquad \qed
\]

{\bf Lemma~4.1}
If $\varphi\in\Phi$ then
\begin{itemize}
\item[(1)] $\liminf\limits_{x\to0^{+}}\;\varphi(x)=\lim\limits_{x\to0^{+}} \; \varphi^{**}(x)$ 
\item[(2)] $\liminf\limits_{x\to+\infty}\;\dfrac{\varphi(x)}{x}=\lim\limits_{x\to+\infty} \; \dfrac{\varphi^{**}(x)}{x}$ 
\end{itemize}

{\it Proof.} Let $\varphi\in\Phi$. Then 
\begin{gather*}
\liminf_{x\to0^{+}} \varphi(x)\ge \liminf_{x\to0^{+}} \varphi^{**}(x)=\lim_{x\to0^{+}} \varphi^{**}(x)\\
\liminf_{x\to+\infty} \frac{\varphi(x)}{x}\ge \liminf_{x\to+\infty} \frac{\varphi^{**}(x)}{x}=\lim_{x\to+\infty} \frac{\varphi^{**}(x)}{x}.
\end{gather*}

Let $a_0\in\R$ and $b_0\in\R$ be such real numbers that $a_0x+b_0\le\varphi(x)$ $\forall x\in(0;+\infty)$.

Thus,
\begin{gather}
\liminf_{x\to0^{+}} \varphi(x)\ge b_0>-\infty\\
\liminf_{x\to+\infty} \frac{\varphi(x)}{x}\ge a_0>-\infty
\end{gather}

Let $b$ be such that $\liminf\limits_{x\to0^{+}} \varphi(x)>b>-\infty$. 

Then we choose a real number $\delta$ in such a way that
$\delta>0$  and 
\[
\inf\limits_{0<x<\delta} \varphi(x)>b.
\]

Therefore 
\begin{align*}
\inf_{x>0}\; \frac{\varphi(x)-b}{x}&\ge\min\Bigl\{ \inf_{0<x<\delta}\; \frac{\varphi(x)-b}{x}\;,\; \inf_{\delta\le x}\; \frac{\varphi(x)-b}{x} \Bigr\}\ge\\
&\ge\min\Bigl\{0,\; \inf_{\delta\le x} \bigl(a_0+\frac{b_0-b}{x} \bigr) \Bigr\}>-\infty.
\end{align*}

Let $a=\min\bigl\{0,\; \inf\limits_{\delta\le x} \bigl(a_0+\frac{b_0-b}{x} \bigr) \bigr\}$.

Thus $(a,b)\in M_{\varphi}$ and consequently $\varphi^{**}(x)\ge ax+b$ $\forall x\in(0;+\infty)$.

So, 
\[
\lim_{x\to0^{+}} \varphi^{**}(x)\ge b
\]
and
\[
\lim_{x\to0^{+}} \varphi^{**}(x)\ge\liminf_{x\to0^{+}} \varphi(x)
\]
because of the choice of the number $b$.

Therefore the assertion (1) of the Lemma~4.1 is proved.

Let $\alpha$ be such that $\liminf\limits_{x\to+\infty} \;\frac{\varphi(x)}{x}>\alpha>-\infty$. 

Then we choose a real number $\Delta$ in such a way that $\Delta>0$ and 
\[
\inf\limits_{x>\Delta} \;\frac{\varphi(x)}{x}>\alpha.
\]

Therfore
\begin{align*}
\inf_{x>0}\; \bigl( \varphi(x)-\alpha x \bigr)&\ge \min\bigl\{ \inf_{0<x\le\Delta}\; ( \varphi(x)-\alpha x ),\; \inf_{x>\Delta}\;( \varphi(x)-\alpha x )\bigr\}\ge\\
&\ge \min\bigl\{ \inf_{0<x\le\Delta}\; ( a_0x+b_0-\alpha x ),\; 0\bigr\}>-\infty
\end{align*}

Let $\beta=\min\bigl\{ \inf\limits_{0<x\le\Delta}\; ( a_0x+b_0-\alpha x ),\; 0\bigr\}$.

Thus, $(\alpha,\beta)\in M_{\varphi}$ and consequently $\varphi^{**}(x)\ge \alpha x+\beta$ $\forall x\in(0;+\infty)$.

So,
\[
\lim_{x\to+\infty} \;\frac{\varphi^{**}(x)}{x}\ge \alpha.
\]
and
\[
\lim_{x\to+\infty} \frac{\varphi^{**}(x)}{x}\ge \liminf_{x\to+\infty} \frac{\varphi(x)}{x}
\]
because of the choice of the number $\alpha$.

Therefore the assertion (2) of the Lemma~4.1 is proved. \qed

{\bf Lemma~4.2}
$\varphi\in\Phi_i$ $\iff$ $\varphi^{**}\in\Phi_i$, ïðè $i=1,2,3$.

{\it Proof.}
The assertion $\varphi\in\Phi_1$ $\iff$ $\varphi^{**}\in\Phi_1$ is proved as (1) of Lemma~4.1.

The proof of the Lemma~4.2 will be completed after proving 
\[
\varphi\in\Phi_3 \iff \varphi^{**}\in\Phi_3.
\]

Let $\varphi\in\Phi_2\cup\Phi_3$ and
\[
\widehat{a}_{\varphi}=\liminf_{x\to+\infty}\;\frac{\varphi(x)}{x}=\lim_{x\to+\infty}\;\frac{\varphi^{**}(x)}{x}.
\]

If $\varphi^{**}\in\Phi_3$ then $\varphi\ge \varphi^{**}$ implies $\varphi\in\Phi_3$.

Now let us suppose that $\varphi\in\Phi_3$. 

Let $a_0$ and $b_0$ be such real numbers that $a_0x+b_0\le\varphi(x)$, $\forall x\in(0;+\infty)$.

Let $b$ be such a real number that $\liminf\limits_{x\to+\infty}\; \bigl( \varphi(x)-\widehat{a}_{\varphi}x\bigr)>b>-\infty$. 

Let the real number $\Delta$ be such that $\Delta>0$ and 
\[
\inf_{x>\Delta}\; \bigl( \varphi(x)-\widehat{a}_{\varphi}x\bigr)>b.
\]

Therefore
\begin{align*}
\inf_{x>0}\; \bigl( \varphi(x)-\widehat{a}_{\varphi}x \bigr)&\ge \min\bigl\{ \inf_{0<x\le\Delta}\; ( \varphi(x)-\widehat{a}_{\varphi}x ),\; \inf_{x>\Delta}\;( \varphi(x)-\widehat{a}_{\varphi}x )\bigr\}\ge\\
&\ge \min\bigl\{ \inf_{0<x\le\Delta}\; ( a_0x+b_0-\widehat{a}_{\varphi}x ),\; b\bigr\}>-\infty
\end{align*}

Let $\widehat{b}=\min\bigl\{ \inf\limits_{0<x\le\Delta}\; ( a_0x+b_0-\widehat{a}_{\varphi}x ),\; b\bigr\}$.

So, $(\widehat{a}_{\varphi},\widehat{b})\in M_{\varphi}$  and consequently 
$\varphi^{**}(x)\ge  \widehat{a}_{\varphi}x+\widehat{b}$, $\forall x\in(0;+\infty)$ and
\[
\liminf_{x\to+\infty} \bigl(\varphi^{**}(x)-\widehat{a}_{\varphi}x\bigr)\ge\widehat{b}>-\infty.
\]

Thus $\varphi^{**}\in\Phi_3$. \qed

{\bf Lemma~4.3} 
Let $\varphi\in\Phi$. If $a$ is such that $a<\widehat{a}_{\varphi}$, then 
\[
\inf\limits_{x>0} \; \bigl( \varphi(x)- ax\bigr)>-\infty
\]
and $\lim\limits_{x\to+\infty} \; \bigl( \varphi(x)- ax\bigr)=+\infty$.

{\it Proof.}
Let $a_0$, $b_0$  be such real numbers that $(a_0,b_0)\in M_{\varphi}$.

Let the real numbers $a$ and $a_1$ be such that
$-\infty<a<a_1<\widehat{a}_{\varphi}$. 

Let the real number $\Delta$ be such that $\Delta>0$ and 
\[
\inf_{x>\Delta} \frac{\varphi(x)}{x}>a_1.
\]

So, $\varphi(x)-ax>(a_1-a)x$, where  $x>\Delta$. 

Therefore $\lim\limits_{x\to+\infty} \; \bigl( \varphi(x)- ax\bigr)=+\infty$ and
\begin{align*}
\inf\limits_{x>0} \; \bigl( \varphi(x)- ax\bigr)&=\min\Bigl\{ \inf\limits_{0<x\le\Delta} \; \bigl( \varphi(x)- ax\bigr);\; \inf\limits_{x>\Delta} \; \bigl( \varphi(x)- ax\bigr)\Bigr\}\ge\\
&\ge \min\Bigl\{ \inf\limits_{0<x\le\Delta} \; \bigl( a_0x+b_0- ax\bigr);\; \inf\limits_{x>\Delta} \; (a_1-a )x \Bigr\}>-\infty.\qquad \qed
\end{align*}

{\bf Lemma~4.4}
Let $\psi:(0;+\infty)\to\R$ be convex on $(0;+\infty)$ and 
\[
\widehat{\psi}(x)=\psi(x)-\psi'(x^{-})x,\quad \forall x\in(0;+\infty),
\]
where $\psi'(x^{-})=\lim\limits_{t\to x^{-}}\frac{\psi(t)-\psi(x)}{t-x}$, $\forall x>0$. 
If $x_1$ and $x_2$ are such that $0<x_1<x_2$ then 
\[
\widehat{\psi}(x_1)\ge\widehat{\psi}(x_2).  
\]

{\it Proof.}
Let $x_1$ and $x_2$ be such real numbers that $0<x_1<x_2$ and let $x_3=\frac{x_1+x_2}{2}$.

Note that 
\begin{itemize}
\item $2\psi(x_3)\le \psi(x_1)+\psi(x_2)$,
\item $f(u,v)=\frac{\psi(u)-\psi(v)}{u-v}$ is such that both functions $f(\cdot,v)$ and $f(u,\cdot)$, each one of them is  a monotone non-decreasing function, where $u>0$, $v>0$, $u\neq v$, and
\[
-\infty<\lim\limits_{t\to x^{-}}\frac{\psi(t)-\psi(x)}{t-x}=\psi'(x^{-})\le\psi'(x^{+})=\lim\limits_{v\to x^{+}}\frac{\psi(v)-\psi(x)}{v-x}<+\infty
\]
where $x>0$.
\end{itemize}

Now, $\widehat{\psi}(x_2)\le\widehat{\psi}(x_1)$ follows from the inequalities  
\begin{align*}
&\widehat{\psi}(x_2)=\psi(x_2)-\psi'(x_2^{-})x_2\le\psi(x_2)-\frac{\psi(x_2)-\psi(x_3)}{x_2-x_3}x_2=\\
&=\bigl(\psi(x_3)-\psi(x_2)\bigr)\frac{2x_2}{x_2-x_1}-\psi(x_2)=\psi(x_3)\frac{2x_2}{x_2-x_1}-\psi(x_2)\frac{x_2+x_1}{x_2-x_1}\le\\
&\le\bigl(\psi(x_2)+\psi(x_1)\bigr)\frac{x_2}{x_2-x_1}-\psi(x_2)\frac{x_2+x_1}{x_2-x_1}=\\
&=\psi(x_1)\frac{x_2+x_1}{x_2-x_1}-\bigl(\psi(x_1)+\psi(x_2) \bigr)\frac{x_1}{x_2-x_1}\le\\
&\le \psi(x_1)\frac{x_2+x_1}{x_2-x_1}-\psi(x_3)\frac{2x_1}{x_2-x_1}=\psi(x_1)-\bigl( \psi(x_3)-\psi(x_1)\bigr)\frac{2x_1}{x_2-x_1}=\\
&=\psi(x_1)-\frac{\psi(x_3)-\psi(x_1)}{x_3-x_1}x_1\le \psi(x_1)-\psi'(x_1^{+})x_1\le \psi(x_1)-\psi'(x_1^{-})x_1=\widehat{\psi}(x_1).\qquad \qed
\end{align*}

{\bf Lemma~4.5}
Let $\psi\in\Phi_2\cup\Phi_3$. If $\psi$ is convex on $(0;+\infty)$ and
\[
\lim_{x\to+\infty}\bigl(\psi(x)-\psi'(x^{-})x\bigr)>-\infty
\]
then $\psi\in\Phi_3$.

{\it Proof.} Note that the limit value exists due to the Lemma~4.4.

Let the real number $\alpha$ be such that  
\[
\lim_{x\to+\infty}\bigl(\psi(x)-\psi'(x^{-})x\bigr)>\alpha>-\infty.
\]

Let thereal number $\Delta$ be such that $\Delta>0$ and $\inf\limits_{x>\Delta} \bigl(\psi(x)-\psi'(x^{-})x \bigr)>\alpha$.

Therefore,
\[
\psi(x)-\frac{\psi(t)-\psi(x)}{t-x}x>\alpha, 
\]
for every $t$ and $x$ such that $\Delta<t<x$. 
 
So,
\[
\frac{\psi(t)-\alpha}{t}\ge \frac{\psi(x)-\alpha}{x}
\]
for every $t$ and $x$ such that $\Delta<t<x$.

Consequently 
\[
\frac{\psi(x)-\alpha}{x}\ge\lim_{x\to+\infty}\frac{\psi(x)-\alpha}{x}=\lim_{x\to+\infty}\frac{\psi(x)}{x}=\widehat{a}_{\psi}
\]
where $x>\Delta$.

Thus $\psi(x)-\widehat{a}_{\psi}x\ge\alpha$, where $x>\Delta$, and 
\[
\lim_{x\to+\infty}\bigl(\psi(x)-\widehat{a}_{\psi}x\bigr)\ge\alpha>-\infty
\]
i.e. $\psi\in\Phi_3$. \qed

{\bf Corollary~4.3} Let $\psi\in\Phi_2$. If $\psi$ is convex on $(0;+\infty)$ then 
\[
\lim_{x\to+\infty}\bigl(\psi(x)-\psi'(x^{-})x\bigr)=-\infty
\]

Note that Corollary~4.3 follows directly from Lemma~4.5.

\section{Proofs of the Main Results}
{\it Proof of Theorem~3.1} Let $\varphi\in\Phi$, $\psi\in\Phi$ and $\psi$ be convex on $(0;+\infty)$.

Note that 
\[
\inf_{x\in(0;+\infty)}\; \bigl(\varphi(x)-\psi(x)\bigr)\ge
\inf_{x\in(0;+\infty)}\; \bigl(\varphi^{**}(x)-\psi(x)\bigr)\eqno(5.1)
\]
because of the inequality $\varphi\ge\varphi^{**}$.

Now we consider two cases 
\begin{itemize}
\item[{\it case 1}] $\inf\limits_{x\in(0;+\infty)}\; \bigl(\varphi(x)-\psi(x)\bigr)=-\infty$. 

So, 
\[
\inf_{x\in(0;+\infty)}\; \bigl(\varphi(x)-\psi(x)\bigr)=
\inf_{x\in(0;+\infty)}\; \bigl(\varphi^{**}(x)-\psi(x)\bigr)=-\infty
\]
\item[{\it case 2}] $c:=\inf\limits_{x\in(0;+\infty)}\; \bigl(\varphi(x)-\psi(x)\bigr)>-\infty$. 

So, 
\[
\varphi(x)\ge\psi(x)+c,\quad \forall x\in(0;+\infty).
\]
and $\psi+c$ is convex minorant of $\varphi$. Therfore
\[
\varphi^{**}(x)\ge\psi(x)+c,\quad \forall x\in(0;+\infty).
\]
i.e. $\inf\limits_{x>0} (\varphi^{**}(x)-\psi(x))\ge c$ and
\[
\inf\limits_{x>0} (\varphi^{**}(x)-\psi(x))\ge \inf\limits_{x\in(0;+\infty)}\; \bigl(\varphi(x)-\psi(x)\bigr)
\]
It follows from here and from inequality~(5.1) that
\[
\inf_{x\in(0;+\infty)}\; \bigl(\varphi(x)-\psi(x)\bigr)=
\inf_{x\in(0;+\infty)}\; \bigl(\varphi^{**}(x)-\psi(x)\bigr).\qquad \qed
\]
\end{itemize}

{\it Proof of Theorem~3.2}
Let the functions $\varphi$ and $\psi$ be such that $\varphi\in\Phi$, $\psi\in\Phi$, $\psi$ is convex on $(0;+\infty)$  and the limit 
\[
\lim\limits_{x\to0^{+}} \bigl(\varphi(x)-\psi(x) \bigr)=+\infty.
\] 

Note that $\lim\limits_{x\to0^{+}}\psi(x)=\psi(0^{+})\in\R\cup\{+\infty\}$ and $\psi(0^{+})>-\infty$.
because $\psi\in\Phi$ and $\psi$ is convex on $(0;+\infty)$.

Therefore, $\varphi(0^{+})=+\infty$ and from Lemma~4.1 we obtain that $\varphi^{**}(0^{+})=+\infty$.

If $\psi(0^{+})<+\infty$ then 
\[
\lim_{x\to0^{+}} (\varphi^{**}(x)-\psi(x))=\varphi^{**}(0^{+})-\psi(0^{+})=+\infty.
\]

So, in order to complete the proof we have to study the alternative when the function $\psi$ is such that $\psi(0^{+})<+\infty$. We will define a new function $\widetilde{\psi}$ that is convex minorant of $\varphi$. 

Let $a_0$ and $b_0$ be such that $(a_0,b_0)\in M_{\varphi}$.

Let $c$ be a real number.

Let $\Delta_1$  be such areal number that $\Delta_1>0$ and
\[
\inf_{0<x<\Delta_1} \bigl(\varphi(x)-\psi(x) \bigr)>c.
\]

Let $\Delta_2$ be such a real number that $\Delta_1>\Delta_2>0$ and
\[
\inf_{0<x<\Delta_2} \bigl( \psi(x)+c-(a_0x+b_0)\bigr)>0.
\]

Let $\Delta_3$ be such a real number that $\Delta_2>\Delta_3>0$ and the function $\psi(x)$ be monotone non-increasing on $(0;\Delta_3)$.

Note that if $x\in(0;\Delta_3)$ then we have the following inequalities with the convex function $\psi$
\[
0\ge \frac{\psi(\Delta_3)-\psi(x)}{\Delta_3}\ge \frac{\psi(\Delta_3)-\psi(x)}{\Delta_3-x}\ge
\limsup_{t\to x^{+}} \; \frac{\psi(t)-\psi(x)}{t-x}=:\psi'(x^{+})
\]
and from $\psi(0^{+})=+\infty$ it follows that
\[
\lim_{x\to0^{+}} \psi'(x^{+})=-\infty.
\] 

Let $\Delta_4$ be such a real number that $\Delta_3>\Delta_4>0$ and
\[
\sup_{0<x<\Delta_4} \psi'(x^{+})<a_0.
\]

Note that if $x\in(0;\Delta_4)$ then
\begin{align*}
&\limsup_{x\to0^{+}} \bigl(\psi'(x^{+})(\Delta_1-x)+\psi(x)+c\bigr) \le\\ 
&\le\limsup_{x\to0^{+}} \bigl(\psi'(x^{+})(0{,}5\Delta_1-x)+\psi(x)+c+\psi'(x^{+})0{,}5\Delta_1\bigr)\le\\
&\le \limsup_{x\to0^{+}}\bigl(\psi(0{,}5\Delta_1)+c+\psi'(x^{+})0{,}5\Delta_1\bigr)=-\infty 
\end{align*}

Let $\Delta_5$ be such a real number that $\Delta_4>\Delta_5>0$ and
\[
\sup_{0<x<\Delta_5} \bigl(\psi'(x^{+})(\Delta_1-x)+\psi(x)+c\bigr)<a_0\Delta_1+b_0.
\] 

Let $x_1$ be such a real number that $x_1\in(0;\Delta_5)$.

Let us set
\begin{align*}
a_1&:=\psi'(x_1^{+}),\\
b_1&:=-\psi'(x_1^{+})x_1+\psi(x_1)+c.
\end{align*}

So,
\[
\psi'(x_1^{+})(x-x_1)+\psi(x_1)+c=a_1x+b_1,\quad \forall x\in(0;+\infty).
\]

Therefore
\begin{align*}
a_1x_1+b_1&=\psi(x_1)+c\ge a_0x_1+b_0,\\
a_1\Delta_1+b_1&< a_0\Delta_1+b_0
\end{align*}

Thus there exists a real number $x_2\in[x_1;\Delta_1]$ such that $a_1x_2+b_1= a_0x_2+b_0$.

And we define the function $\widetilde{\psi}$
\[
\widetilde{\psi}(x):=\begin{cases}
\psi(x)+c,&\quad x\in(0;x_1)\\
a_1x+b_1,&\quad x\in[x_1;x_2]\\
a_0x+b_0,&\quad x\in(x_2;+\infty).
\end{cases}
\]

The function $\widetilde{\psi}$ is convex on $(0;+\infty)$ because it is continuous, $\psi'(x_1^{-})\le a_1\le a_0$ and $\psi+c$ is convex on $(0;x_1)$. 

Furthermore
\begin{align*}
\widetilde{\psi}(x)&=\psi(x)+c\le\varphi(x), \quad \forall x\in(0;x_1)\\
\widetilde{\psi}(x)&=a_1x+b_1\le \psi(x)+c\le\varphi(x), \quad \forall x\in[x_1;x_2]\\
\widetilde{\psi}(x)&=a_0x+b_0\le\varphi(x),  \quad \forall x\in(x_2;+\infty)
\end{align*}
i.e. $\widetilde{\psi}$ is convex minorant of $\varphi$.

Therfore $\varphi^{**}(x)\ge \widetilde{\psi}(x)$, $\forall x\in(0;+\infty)$.

Thus $\varphi^{**}(x)\ge \psi(x)+c$ $\forall x\in(0;x_1)$ and
\[
\liminf_{x\to0^{+}} \bigl(\varphi^{**}(x)- \psi(x)\bigr)\ge c.
\]
which, accordingly to the choice of the number $c$, implies that
\[
\lim_{x\to0^{+}} \bigl(\varphi^{**}(x)- \psi(x)\bigr)=+\infty.\qquad\qed
\]

{\it Proof of Theorem~3.3}
Let the functions $\varphi$ and $\psi$ be such that $\varphi\in\Phi$, $\psi\in\Phi\setminus\Phi_3$, 
$\psi$ be convex on $(0;+\infty)$ and the limit value $\lim\limits_{x\to+\infty} \bigl(\varphi(x)-\psi(x) \bigr)=+\infty$.

Note that by Lemma~4.1
\begin{align*}
\widehat{a}_{\varphi}&=\liminf_{x\to+\infty}\; \frac{\varphi(x)}{x}=\lim_{x\to+\infty}\; \frac{\varphi^{**}(x)}{x},\\
\widehat{a}_{\psi}&=\liminf_{x\to+\infty}\; \frac{\psi(x)}{x}=\lim_{x\to+\infty}\; \frac{\psi^{**}(x)}{x}
\end{align*}

Let the real number $\Delta$ be such that $\Delta>0$ and 
\[
\inf_{x>\Delta}\bigl(\varphi(x)-\psi(x)\bigr)>0.
\]

Therefore $\varphi(x)\ge\psi(x)$,  $\forall x\in(\Delta,+\infty)$, and
\[
\widehat{a}_{\varphi}=\liminf_{x\to+\infty}\; \frac{\varphi(x)}{x}\ge \liminf_{x\to+\infty}\; \frac{\psi(x)}{x} =\widehat{a}_{\psi}.
\]

In order to prove Theorem~3.3 we consider the following four cases 
\begin{itemize}
\item[{\it case 1. }] $\widehat{a}_{\varphi}>\widehat{a}_{\psi}$
\item[{\it case 2. }] $\widehat{a}_{\varphi}=\widehat{a}_{\psi}$
  \begin{itemize}
  \item[{\it case 2.1. }] $\varphi\in\Phi_3$
  \item[{\it case 2.2. }] $\varphi\in\Phi_2$
  \item[{\it case 2.3. }] $\varphi\in\Phi_1$.
  \end{itemize}  
\end{itemize}

{\it Case 1.} Let the functions $\varphi$ and $\psi$ be such that $\widehat{a}_{\varphi}>\widehat{a}_{\psi}$.

Let the real numbers $a_1$ and $a_2$ be such that
\[
\widehat{a}_{\varphi}>a_1>a_2>\widehat{a}_{\psi}.
\]

Let the real number $\Delta_1$ be such that $\Delta<\Delta_1$ and
\[
\frac{\varphi^{**}(x)}{x}>a_1>a_2>\frac{\psi(x)}{x},\quad \forall x\in(\Delta_1,+\infty).
\]

So, $\varphi^{**}(x)-\psi(x)>(a_1-a_2)x$, $\forall x\in(\Delta_1,+\infty)$, and
\[
\lim_{x\to+\infty}\bigl(\varphi^{**}(x)-\psi(x)\bigr)=+\infty.
\]

Thus the case 1 is closed.

{\it Case 2.1.} Let the functions $\varphi$ and $\psi$ be such that  $\widehat{a}_{\varphi}=\widehat{a}_{\psi}$ and $\varphi\in\Phi_3$.

So,
\begin{itemize}
\item Lemma~3.2 implies that $\varphi^{**}\in\Phi_3$. 
\item $\psi\in\Phi_2$ and
\[
\liminf_{x\to+\infty}\bigl(\psi(x)-\widehat{a}_{\psi}x\bigr)
=\lim_{x\to+\infty}\bigl(\psi(x)-\widehat{a}_{\psi}x\bigr)=-\infty
\]
because of the convexity of $\psi$.
\end{itemize}

We claim that 
\[
\inf\limits_{x>0}\; \bigl(\varphi^{**}(x)-\widehat{a}_{\varphi}x\bigr)>-\infty.\eqno(5.2)
\]

Indeed, let us choose the real numbers $b$, $\Delta_2$, $a_0$ and $b_0$ 
\begin{itemize}
\item the number $b$ is such that 
\[
\liminf_{x\to+\infty}\bigl(\varphi^{**}(x)-\widehat{a}_{\varphi}x\bigr)>b
\] 
\item the number $\Delta_2$ is such that $\Delta_2>0$ and
\[
\inf_{x>\Delta_2}\bigl(\varphi^{**}(x)-\widehat{a}_{\varphi}x\bigr)>b
\]
\item the numbers $a_0$ and $b_0$ are such that $(a_0,b_0)\in M_{\varphi}$ and therefore 
\[
a_0x+b_0\le\varphi^{**}(x)\quad\forall x\in(0,+\infty).
\]
\end{itemize}

So
\begin{align*}
\inf\limits_{x>0}\; \bigl(\varphi^{**}(x)-\widehat{a}_{\varphi}x\bigr)&=
\min\Bigl\{\inf\limits_{0<x\le\Delta_2}\; \bigl(\varphi^{**}(x)-\widehat{a}_{\varphi}x\bigr);\;
\inf\limits_{x>\Delta_2}\; \bigl(\varphi^{**}(x)-\widehat{a}_{\varphi}x\bigr) \Bigr\}\ge\\
&\ge \min\Bigl\{\inf\limits_{0<x\le\Delta_2}\; \bigl(a_0x+b_0-\widehat{a}_{\varphi}x\bigr);\;
b \Bigr\}>-\infty
\end{align*}
and the claim~(5.2) is proved.

Thus,
\begin{align*}
\liminf_{x\to+\infty}\bigl(\varphi^{**}(x)-\psi(x)\bigr)&=
\liminf_{x\to+\infty}\Bigl((\varphi^{**}(x)-\widehat{a}_{\varphi}x)+(\widehat{a}_{\psi}x-\psi(x))\Bigr)\ge\\
&\ge \inf_{x>0}(\varphi^{**}(x)-\widehat{a}_{\varphi}x) + \liminf_{x\to+\infty}(\widehat{a}_{\psi}x-\psi(x))=\\
&=\inf_{x>0}(\varphi^{**}(x)-\widehat{a}_{\varphi}x) + \lim_{x\to+\infty}(\widehat{a}_{\psi}x-\psi(x))=+\infty
\end{align*}

The case 2.1 is closed.

{\it Case 2.2 and Case 2.3} Let $c$ be a real number.

Let the number $\Delta$ be such that $\Delta>0$ and $\inf\limits_{x>\Delta}(\varphi(x)-\psi(x))>c$.

Let the numbers $a_0$ and $b_0$ be such that $(a_0,b_0)\in M_{\varphi}$.

In the present cases, the assumptions imply that $\varphi\in\Phi_2$ and $+\infty>\widehat{a}_{\varphi}>a_0$.

So, $\widehat{a}_{\psi}=\widehat{a}_{\varphi}>a_0$ and by Lemma~4.3  
\[
\lim_{x\to+\infty}(\psi(x)-a_0x)=+\infty.
\]

Let the number $\Delta_1$ be such that $\Delta_1>\Delta$ and 
\[
\inf_{x>\Delta_1}\bigl( \psi(x)+c-(a_0x+b_0)\bigr)>0.
\]

Note that the convex function $\psi$ belongs to $\Phi_2$. Therefore 
\[
\psi'(x^{-})\le\psi'(x^{+})<\widehat{a}_{\psi},\quad \forall x>0.
\]

Furthermore if we fix an $x'$ and let $x$ is such that $+\infty>x>x'>0$ then 
\begin{align*}
\frac{\psi(x')-\psi(x)}{x'-x}&\le \lim_{t\to x^{-}} \frac{\psi(t)-\psi(x)}{t-x}=\psi'(x^{-}) \le\widehat{a}_{\psi}\\
\implies \lim_{x\to+\infty}\frac{\frac{\psi(x)}{x}-\frac{\psi(x')}{x}}{1-\frac{x'}{x}}&=\lim_{x\to+\infty}\frac{\psi(x')-\psi(x)}{x'-x}\le \lim_{x\to+\infty}\psi'(x^{-}) \le\widehat{a}_{\psi}
\end{align*}
and thus
\[
\lim_{x\to+\infty}\psi'(x^{-}) =\widehat{a}_{\psi}\eqno(5.3)
\]

Let the real number $\Delta_2$ be such that $\Delta_2>\Delta_1$ and 
\[
\inf\limits_{x>\Delta_2} \psi'(x^{-})>a_0.
\]

We claim that there exists $\Delta_3$ such that $\Delta_3>\Delta_2$ and 
\[
\psi'(x^{-})(\Delta-x)+\psi(x)+c<a_0\Delta+b_0,\quad \forall x>\Delta_3.
\]

Argument for this claim in the case 2.2 is different when compared to the corresponding one in the case 2.3.
\begin{itemize}
\item[{\it case 2.2}] in this case $\widehat{a}_{\psi}<+\infty$ and Corollary~4.3 applied to the function $\psi$ imply that 
\[
\lim\limits_{x\to+\infty}(\psi(x)-\psi'(x^{-})x)=-\infty
\]
Therefore by (5.3) we obtain
\[
\lim_{x\to+\infty}\bigl(\psi'(x^{-})(\Delta-x)+\psi(x)+c\bigr)=-\infty
\]
\item[{\it case 2.3}] now, $\lim\limits_{x\to+\infty} \psi'(x^{-})=\widehat{a}_{\psi}=+\infty$
and
\begin{align*}
\psi'(x^{-})(\Delta-x)+\psi(x)+c&\le \psi'(x^{-})(2\Delta-x)+\psi(x)+c-\psi'(x^{-})\Delta\le\\
&\le \psi(2\Delta) +c-\psi'(x^{-})\Delta
\end{align*}
because of $\psi'(x^{-})(t-x)+\psi(x)\le\psi(t)$ with  $x>0$ and $t>0$.

So,
\[
\lim_{x\to+\infty}(\psi'(x^{-})(\Delta-x)+\psi(x)+c)=-\infty
\]
\end{itemize}

Thus the claim is proved and let the number $\Delta_3$ is such that $\Delta_3>\Delta_2$ and 
\[
\psi'(x^{-})(\Delta-x)+\psi(x)+c<a_0\Delta+b_0,\quad \forall x>\Delta_3.
\].

Let the real number $x_1$ be such that $x_1>\Delta_3$ and let us set 
\begin{align*}
a_1&=\psi'(x_1^{-}),\\
b_1&=-\psi'(x_1^{-})x_1+\psi(x_1)+c.
\end{align*}
Note that $a_1>a_0$. 

Thus
\begin{align*}
&a_1x+b_1\le\psi(x)+c,	\quad \forall x\in(0;+\infty)\\
&a_1x_1+b_1=\psi(x_1)+c\ge a_0x_1+b_0,\\
&a_1\Delta+b_1<a_0\Delta+b_0.
\end{align*}

Let the real number $x_2$ be such that $x_2\in(\Delta;x_1]$ and
\[
a_1x_2+b_1=a_0x_2+b_0.
\]

We define the function $\widetilde{\psi}:(0;+\infty)\to\R$
\[
\widetilde{\psi}(x)=\begin{cases}
a_0x+b_0,&\quad x\in(0;x_2]\\
a_1x+b_1,&\quad x\in(x_2;x_1]\\
\psi(x)+c,&\quad x\in(x_1;+\infty).
\end{cases}
\]

The function $\widetilde{\psi}$ is convex on $(0;+\infty)$ because it is continuous, $a_0\le a_1\le\psi'(x_1^{-})$ and $\psi+c$ is convex on $(x_1;+\infty)$.

Note that $\widetilde{\psi}(x)\le\varphi(x)$, $\forall x\in(0;+\infty)$.

Therfore $\widetilde{\psi}(x)\le\varphi^{**}(x)$, $\forall x\in(0;+\infty)$.

Thus with $x>x_1$ we have $\psi(x)+c\le \varphi^{**}(x)$ and
\[
\liminf_{x\to+\infty} \bigl(\varphi^{**}(x)-\psi(x)\bigr)\ge c
\]

Accordingly to the choice of the number $c$
\[
\lim_{x\to+\infty} \bigl(\varphi^{**}(x)-\psi(x)\bigr)=+\infty.\qquad\qed
\]

\section{Application}
In this section we apply the theorems 3.1, 3.2, 3.3 to the theory of spaces $H_{v}(G)$ and $H_{{v}_0}(G)$.

We make use of the following notation
\begin{align*}
Mf(y)&=\sup_{x\in(-\infty;+\infty)}\abs{f(x+iy)},\\
\psi_f(y)&=\ln Mf(y),\qquad \forall y>0,f\in \Lambda(p),
\end{align*}
where $f$ is a holomorphic function defined on the upper half plane $G$.

Note that
\[
(-1)\ln \parallel f\parallel_v=\inf_{y>0}\bigl(\varphi_v(y)-\psi_f(y) \bigr)
\]

Here we reformulate our results from the e-preprint [6]

{\bf Theorem A.}[6, Th.~1.2] If $\varphi$ meets the condition~(1.1') then 
\[
H_v(G)\ne\{0\} \iff \varphi\in\Phi.
\]
where $v=e^{(-1)\varphi}$.

{\bf Theorem B.}[6, Th.~1.3] If $\varphi$ meets the condition~(1.1') then 
\[
H_{v_0}(G)\ne\{0\} \iff \left|\begin{array}{l}
\varphi\in\Phi,\\
\varphi(0^{+})=+\infty
\end{array} \right.
\]
where $v=e^{(-1)\varphi}$.

{\bf Theorem C}[6, Th.~1.4] If $\varphi$ meets the condition~(1.1') and $H_{v_0}(G)\ne\{0\}$ then
\[
\psi_f\in\Phi\setminus\Phi_3\qquad \forall f\in H_{v_0}(G)\setminus \{0\}
\]
where $v=e^{(-1)\varphi}$.

Note that $\psi_f$ is convex on $(0;+\infty)$ and $\psi_f\in\Phi$, $\forall f\in H_{v}(G)\setminus \{0\}$. 

Now in this section we prove two new theorems.

{\bf Theorem~6.1} If $\varphi$ meets the condition~(1.1') and $\varphi\in\Phi$ then
\[
\bigl(H_{v}(G),\parallel \cdot\parallel_v\bigr)\equiv\bigl(H_{w}(G),\parallel \cdot\parallel_w\bigr)
\]
where $v=e^{(-1)\varphi}$ and $w=e^{(-1)\varphi^{**}}$. 

{\it Proof.} Let $v=e^{(-1)\varphi}$ and $w=e^{(-1)\varphi^{**}}$.

Note that
\begin{itemize} 
\item $\varphi>\varphi^{**}$ $\implies$ $\varphi^{**}$ meets the condition~(1.1') and 
\item  $\varphi\in\Phi$ $\implies$ $\varphi^{**}\in\Phi$, because $M_{\varphi^{**}}=M_{\varphi}\neq\emptyset$. 
\end{itemize}
Thus, $H_v(G)\ne\{0\}$ and $H_w(G)\ne\{0\}$  by Theorem~A. 
 
$H_v(G)\supset H_w(G)$, because $\parallel f\parallel_v\le\parallel f\parallel_w<+\infty$, $\forall f\in H_w(G)$. 

Note that $\forall f\in H_v(G)\ne\{0\}$ the function $\psi_f=\ln Mf$ is convex on $(0;+\infty)$ and $\psi_f\in\Phi$.
Therefore by Theorem~3.1
\[
\inf_{x\in(0;+\infty)}\; \bigl(\varphi(x)-\psi_f(x)\bigr)=
\inf_{x\in(0;+\infty)}\; \bigl(\varphi^{**}(x)-\psi_f(x)\bigr)
\]

Thus $f\in H_w(G)\ne\{0\}$ and $\parallel f\parallel_v=\parallel f\parallel_w$. \qed

{\bf Theorem~6.2} If $\varphi$  meets the condition~(1.1'), $\varphi\in\Phi$ and $\varphi(0^{+})=+\infty$ then 
\[
\bigl(H_{v_0}(G),\parallel \cdot\parallel_v\bigr)\equiv\bigl(H_{w_0}(G),\parallel \cdot\parallel_w\bigr)
\]
where $v=e^{(-1)\varphi}$ and $w=e^{(-1)\varphi^{**}}$. 

{\it Proof.}  Let $v=e^{(-1)\varphi}$ and $w=e^{(-1)\varphi^{**}}$.

Note that
\begin{itemize} 
\item $\varphi>\varphi^{**}$ $\implies$ $\varphi^{**}$ meets the condition~(1.1') and 
\item  $\varphi\in\Phi$ $\implies$ $\varphi^{**}\in\Phi$, because $M_{\varphi^{**}}=M_{\varphi}\neq\emptyset$
\item $\varphi^{**}(0^{+})=\varphi(0^{+})=+\infty$, because of Lemma~4.1(1).
\end{itemize}
Thus, $H_{v_0}(G)\ne\{0\}$ and $H_{w_0}(G)\ne\{0\}$  by Theorem~B. 

$H_{v_0}(G)\supset H_{w_0}(G)$, because 
\[
0\le v(iy)\abs{f(x+iy)}\le w(iy)\abs{f(x+iy)},
\]
$\forall f\in H_{w_0}(G)$ and $\forall x\in(-\infty;+\infty)$, $\forall y\in(0;+\infty)$. 

Note that $\parallel f\parallel_v=\parallel f\parallel_w$, $\forall f\in H_{v_0}(G)\ne\{0\}$ because of Theorem~6.1. 

Now we have to prove that $f\in H_{w_0}(G)\ne\{0\}$, $\forall f\in H_{v_0}(G)\ne\{0\}$.

Let $f\in H_{v_0}(G)\ne\{0\}$. Accordingly to the definition of $H_{v_0}(G)$
\[
\lim_{\mathcal K \uparrow G} \sup_{z\in G\setminus \mathcal K} v(z)\abs{f(z)}=0
\]
where $\mathcal K\subset G$ and $\mathcal K$ is compact.
So,
\[
\lim_{y\to0^{+}} v(iy)Mf(y)=0,\quad 
\lim_{y\to+\infty} v(iy)Mf(y)=0
\]
and after reformulation
\begin{align*}
&\lim_{y\to0^{+}} \bigl(\varphi(y)-\psi_f(y)\bigr)=+\infty,\\
&\lim_{y\to+\infty} \bigl(\varphi(y)-\psi_f(y)\bigr)=+\infty
\end{align*}

Note that $\psi_f\in\Phi\setminus\Phi_3$ by Theorem~C.

Therefore by Theorem~3.2 and Theorem~3.3
\begin{align*}
&\lim_{y\to0^{+}} \bigl(\varphi^{**}(y)-\psi_f(y)\bigr)=+\infty,\\
&\lim_{y\to+\infty} \bigl(\varphi^{**}(y)-\psi_f(y)\bigr)=+\infty
\end{align*}
i.e.
\[
\lim_{y\to0^{+}} w(iy)Mf(y)=0,\quad
\lim_{y\to+\infty} w(iy)Mf(y)=0
\]

Let the real number $\varepsilon$ be such that $\varepsilon>0$.

Let the real number $c$ be such that $c>1$ and
\[
\sup_{y<\frac1c} w(iy)Mf(y)<\varepsilon,\quad
\sup_{y>c} w(iy)Mf(y)<\varepsilon
\]

Let 
\[
m=\frac{\sup\limits_{\frac1c\le y\le c}w(iy)}{\inf\limits_{\frac1c\le y \le c}v(iy)}
\]
$m<+\infty$ because $\varphi^{**}\in\Phi$ and therefore $\inf_{\frac1c\le x\le c}\varphi^{**}(x)>-\infty$, $\forall c>1$.

So, $0<m<+\infty$.

Accordingly to the definition of $H_{v_0}(G)$
there exist two numbers $x_1$ and $c_1$ such that $x_1>0$, $c_1>c$ 
and the compact 
\[
\mathcal K_1=\{x+iy\,|\; -x_1\le x\le x_1, \frac1{c_1}\le y\le c_1\}
\]
meets the condition
\[
\sup_{x+iy\in G\setminus \mathcal K_1} v(iy)\abs{f(x+iy)}\le\frac{\varepsilon}{m}
\]  

Let $\mathcal K=\{x+iy\,|\; -x_1\le x\le x_1, \frac1c\le y\le c\}$

Thus
\begin{align*}
\sup_{x+iy\in G\setminus \mathcal K} w(iy)\abs{f(x+iy)}&=
   \max\Bigl\{\;
         \sup_{y<\frac1c} w(iy)Mf(y), \;
         \sup_{\substack{\abs{x}>x_1,\\ \frac1c\le y\le c }}w(iy)Mf(y) ,\;
         \sup_{y>c} w(iy)Mf(y)\;
       \Bigr\}\le\\
   &\le\max\Bigl\{\;
         \varepsilon, \;
         \sup_{\substack{\abs{x}>x_1,\\ \frac1c\le y\le c }}v(iy)m\abs{f(x+iy)},\;
         \varepsilon \;
       \Bigr\}\le\varepsilon
\end{align*}
and therefore $f\in H_{w_0}(G)$. \qed


\section{REFERENCES}

\small
\begin{enumerate} \frenchspacing

\item{ \vskip-2pt
M. Ardalani and W. Lusky,
    Weighted Spaces of Holomorphic Functions
    on the Upper Halfplane,
    Math. Scand. 111 (2012), 244–260,
    Zbl~1267.30111
}

\item{ \vskip-2pt
A.~Harutyunyan and W.~Lusky,
    A remark on the isomorphic classification of weighted spaces of holomorphic functions on the upper half plane,
    Ann. Univ. Sci. Budap. Sect. Comp.,
    39,
    125--135,
    2013,
    Zbl~1289.46045
}

\item{ \vskip-2pt
J.~Bonet, P.~Domanski, M.~Lindstrom and J.~Taskinen,
    Composition Operators Between Weighted Banach
    Spaces of Analytic Functions,
    J.~Austral. Math. Soc. (Series A) 64 (1998), 101-118,
    Zbl~0912.47014
}

\item{ \vskip-2pt
M.~Contreras and A.~Hernandez-Diaz,
   Weighted Composition Operators in Weighted
   Banach Spaces of Analytic Functions,
   J.~Austral. Math. Soc. (Series A) 69 (2000), 41-60
   Zbl~0990.47018
}

\item{\vskip-2pt K.~Bierstedt, J.~Bonet, J. Taskinen,
   Associated Weights and Spaces of Holomorphic Functions,
   Stud. Math. 127 (2) (1998), 137--168,
   Zbl~0934.46027
}

\item{\vskip-2pt
 Martin At. Stanev, Weighted Banach spaces of holomorphic
functions in the upper half plane,  (1999), e-preprint http://arxiv.org/abs/math/9911082
}

\item{\vskip-2pt
Martin Stanev,
Log-convexity of the Weight of a Weighted Function Space
"Complex Analysis and Applications '13"
(International Memorial Conference for the 100th Anniversary of Acad. Ljubomir Iliev)
IMI - Sofia, 31 Oct. - 2 Nov. 2013,
http://www.math.bas.bg/complan/caa13/
}

\item{\vskip-2pt
Martin Stanev,
Weighted Banach Spaces of Holomorphic Functions with
Log-Concave Weight Function,
Jubilee Conference
125 years of Mathematics and Natural Sciences at Sofia University "St. Kliment Ohridski",
December  6-7, 2014,
http://125years.fmi.uni-sofia.bg
}

\end{enumerate}
%
\vskip20pt
 \footnotesize
\begin{flushleft}
Martin At. Stanev \\
Dep. of Mathematics and Physics \\
University of Forestry \\
10 blvd.~K.~Ohridski, BG-1756 Sofia\\
BULGARIA \\
e-mail: martin\_stanev@yahoo.com
\end{flushleft}

\end{document}